\documentclass{amsart}[11pt,leqno] \usepackage{amssymb,enumerate,amscd,tikz,url}

\theoremstyle{plain}
\newtheorem{theo}[equation]{Theorem}
\newtheorem{lem}[equation]{Lemma}
\newtheorem{cor}[equation]{Corollary}
\newtheorem{prop}[equation]{Proposition}
\theoremstyle{definition}
\newtheorem{Def}[equation]{Definition}
\newtheorem{Ex}[equation]{Example}
\newtheorem{Rem}[equation]{Remark} 
\newtheorem{Not}[equation]{Notation}

\newcommand{\F}{{\mathbb F}}

\newcommand{\Q}{{\mathbb Q}}

\newcommand{\T}{{\mathbb T}}
\newcommand{\W}{{\mathbb W}}
\newcommand{\Z}{{\mathbb Z}}

\newcommand{\bfa}{{\mathbf a}}

\newcommand{\bfs}{{\bf s}}

\newcommand{\Mu}{\boldsymbol \mu}
\newcommand{\fraksp}{{\mathfrak{sp}}}

\newcommand{\gA}{{\mathfrak A}}
\newcommand{\gE}{{\mathfrak E}}
\newcommand{\gO}{{\mathfrak O}}

\newcommand{\gf}{{\mathfrak f}}

\newcommand{\gM}{{\mathfrak M}}

\newcommand{\gn}{{\mathfrak n}}

\newcommand{\gp}{{\mathfrak p}}

\newcommand{\gR}{{\mathfrak R}}
\newcommand{\gS}{{\mathfrak S}}

\newcommand{\cB}{{\mathcal B}}

\newcommand{\cE}{{\mathcal E}}
\newcommand{\cF}{{\mathcal F}}
\newcommand{\cG}{{\mathcal G}}
\newcommand{\cH}{{\mathcal H}}

\newcommand{\cL}{{\mathcal L}}
\newcommand{\cM}{{\mathcal M}}

\newcommand{\cO}{{\mathcal O}}

\newcommand{\cS}{{\mathcal S}}

\newcommand{\cV}{{\mathcal V}}
\newcommand{\cW}{{\mathcal W}}

\newcommand{\rF}{{\rm F}}
\newcommand{\rM}{{\rm M}}
\newcommand{\rL}{{\rm L}}
\newcommand{\rV}{{\rm V}}

\newcommand{\CW}{\operatorname{CW}}

\newcommand{\Gal}{\operatorname{Gal}}
\newcommand{\Hom}{\operatorname{Hom}}

\newcommand{\Ext}{\operatorname{Ext}}
\newcommand{\GL}{{\rm GL}}
\newcommand{\Mat}{{\rm Mat}}
\newcommand{\SL}{{\rm SL}}
\newcommand{\GS}{{\rm GSp}}
\newcommand{\SP}{{\rm Sp}}
\newcommand{\spn}{{\rm span}}

\newcommand{\Sym}{\operatorname{Sym}}

\newcommand{\ord}{\operatorname{ord}}

\newcommand{\0}{\vec{0}}

\newcommand{\bino}[2]{\left( \begin{smallmatrix} {#1} \\ {#2} \end{smallmatrix} \right)}

\newcommand{\ov}[1]{\overline{#1}}

\newcommand{\fdeg}[2]{[{#1}\!:\!{#2}]}

\newcommand{\lr}[1]{\langle{#1}\rangle}

\makeatletter
\renewcommand*\l@subsection{\@tocline{2}{0pt}{30pt}{0pt}{}}
\makeatother

\theoremstyle{plain}

\begin{document}
\title[Large Galois image and non-existence]{Large 2-adic Galois image and non-existence of certain abelian surfaces over $\Q$}

\begin{abstract} Motivated  by our arithmetic applications, we required some tools that might be of independent interest. 

Let $\cE$ be  an absolutely irreducible group scheme of rank $p^4$  over $\Z_p$.  We provide a complete description of the Honda systems of $p$-divisible groups $\cG$ such that $\cG[p^{n+1}]/\cG[p^n] \simeq \cE$ for all $n$.  Then we find a bound for the abelian conductor of the second layer $\Q_p(\cG[p^2])/\Q_p(\cG[p])$, stronger in our case than can be deduced from Fontaine's bound. 

Let $\pi\!: \, \SP_{2g}(\Z_p) \to \SP_{2g}(\F_p)$ be the reduction map and let $G$ be a closed subgroup of $\SP_{2g}(\Z_p)$ with $\ov{G} = \pi(G)$ irreducible and generated by transvections.  We fill a gap in the literature by showing that if $p=2$ and $G$ contains a transvection, then $G$ is as large as possible in $\SP_{2g}(\Z_p)$ with given reduction $\ov{G}$, i.e. $G = \pi^{-1}(\ov{G})$.  

One simple application arises when $A = J(C)$ is the Jacobian of a hyperelliptic curve $C\!: \, y^2 + Q(x)y = P(x)$, where $Q(x)^2 + 4P(x)$ is irreducible in $\Z[x]$ of degree $m=2g+1$ or $2g+2$, with Galois group $\cS_m \subset \SP_{2g}(\F_2)$.  If the Igusa discriminant $I_{10}$ of $C$ is odd and some prime $q$ exactly divides $I_{10}$, then $G = \Gal(\Q(A[2^\infty])/\Q)$ is $\tilde{\pi}^{-1}(\cS_m)$, where $\tilde{\pi}\!: \, \GS_{2g}(\Z_p) \to \SP_{2g}(\F_p)$.  

When $m = 5$, $Q(x) = 1$ and $I_{10} = N$ is a prime, $A = J(C)$ is an example of a {\em favorable} abelian surface.  We use the machinery above to obtain non-existence results for certain favorable abelian surfaces, even for large $N$.
\end{abstract}

\author[A. Brumer]{Armand Brumer}
\address{Department of Mathematics, Fordham University, Bronx, NY 10458}
\email{brumer@fordham.edu}
\author[K. Kramer]{Kenneth Kramer}
\address{Department of Mathematics, Queens College (CUNY), Flushing, NY 11367;  Department of Mathematics, The Graduate Center of CUNY, New York, NY 10016}
\email{kkramer@qc.cuny.edu}

\maketitle 

\tableofcontents 

\numberwithin{equation}{section}
\section{Introduction} 
Let $A_{/\Q}$ be a $g$-dimensional abelian variety and  $G = \Gal(\Q(A[p^\infty])/\Q)$  the Galois group of its $p$-division tower. Serre's work on the open image problem for abelian varieties has stimulated a large literature. For instance,  \cite{ALS, Hall, KM, Vas, ZS} show that, under suitable hypotheses, $G$ is an open subgroup of $\GS_{2g}(\Z_p),$ at least for large $p.$  We concentrate  on $p=2,$  for which more residual images exist.  Theorem \ref{fullpre} describes our group theoretical conclusions.  As usual, suitable abelian varieties thereby give rise to large Galois extensions with controlled ramification, as in Proposition \ref{BigGal}.  

Given an integer $N$ and a group scheme $\cE$ over $\Z[\frac{1}{N}]$ of exponent $p,$ one may ask for the existence (or  even the  uniqueness  up to isogeny) of an abelian variety $A$ with $A[p]\simeq \cE$.  In \cite{BK1}, we found non-existence criteria when $\cE$ is reducible.  In this paper, we treat non-existence criteria when $\dim A = 2$ and $\cE_{\vert \Z_p}$ is absolutely irreducible.  Then $A$ has a polarization of degree prime to $p$ and $\Gal(\Q(A[p^\infty])/\Q)$ is contained in $\GS_4(\Z_p).$ This requires a delicate study of the extensions $\cW$ of $\cE$ by $\cE$ of exponent $p^2.$  Non-existence of  $\Q(\cW)/\Q$ implies that of $A$ with $A[p^2] \simeq \cW.$ 

Let $\cE$ be an absolutely irreducible group scheme of rank $p^4$  over $\Z_p$.  In \S\ref{pdiv}, we give a complete description of the Honda systems of $p$-divisible groups $\cG$ such that $\cG[p^{n+1}]/\cG[p^n] \simeq \cE$ for all $n$.   In \S\ref{expp2}, we study the field of points of $\cW= \cG[p^2]$, thereby obtaining a bound for the abelian conductor of $\Q_p(\cW)/\Q_p(\cE),$ stronger in our special case than can be deduced from Fontaine's bound.   When $A$ is the Jacobian of a genus 2 curve over $\Q_2$, the parameters associated to $A[4]$ are determined in Proposition \ref{FindParams}.

For the global applications in \S \ref{global}, let $p=2$ and recall the following definition.

\begin{Def}[\cite{BK3}]  \label{fav}
A quintic field is {\em favorable} if its discriminant is $\pm 16N$ with $N$ prime  and the ramification index over the prime 2 is 5.  An abelian surface $A_{/\Q}$ of conductor $N$ is {\em favorable} if its 2-division field is the Galois closure of a favorable quintic field.
\end{Def}

If $A$ is favorable, then the image in $\SP_4(\F_2)$ of the representation of $G_\Q$ on $A[2]$ is $O_4^-(\F_2) \simeq \cS_5$, with transvections corresponding to transpositions.  In addition, $A[2]$ is biconnected, absolutely simple and Cartier self-dual over $\Z_2$.  Let $\gS =\pi^{-1}(\cS_5)$, where $\pi\!: \, \GS_4(\Z/4\Z) \to \SP_4(\F_2)$ is the natural projection.  By Remark \ref{favS}, $\Q(A[4])$ is a {\em favorable} $\gS$-field with $F = \Q(A[2])$, as in the following definition. 

\begin{Def} \label{Sfield}
Fix the Galois closure $F$ of a favorable quintic field of discriminant $\pm 16N$ and let $L$ be a field containing $F$.  Then $L$ is a {\em favorable} $\gS$-field if $L/\Q$ is Galois, with $\Gal(L/\Q) \simeq \gS$ and
\begin{enumerate}[i)]
\item $L/\Q$ is unramified outside $2 N \infty$;  \vspace{2 pt}
\item the abelian conductor exponent at primes over 2 in $L/F$ is 6 and \vspace{2 pt}
\item the inertia group at each prime over $N$ in $L/\Q$ is generated by a transvection.
\end{enumerate}
\end{Def} 

Proposition \ref{STest} gives a testable ray class field criterion necessary for the existence of a favorable $\gS$-field $L$.   This explains the non-existence results in \cite[Table 3]{BK3}. 

\section{A large image result} \label{FULL}   
\numberwithin{equation}{subsection}

\subsection{Review}  \label{FullReview}
For closed subgroups $\Gamma$ of $\GL_m(\Z_p)$, set $\Gamma^{(n)} = \{ g \in \Gamma \, \vert \, g \equiv I_m \, (p^n) \}$ and $\ov{\Gamma} = \Gamma/\Gamma^{(1)}$.  A closed subgroup $G$ of $\Gamma$ is {\em saturated in} $\Gamma$ if $G^{(1)} = \Gamma^{(1)}$, so that $G$ is as large as possible in $\Gamma$, subject to its reduction being $\ov{G}$.  When there is a symplectic pairing $[ \hspace{3 pt}, \hspace{3 pt}]\!: \, M \times M \to \Z_p$ on $M =\Z_p^{2g}$, {\em transvections} in $\SP(M)$ have the form $\sigma(x) = x - \lambda \, [y,x] \, y$,  with $y$ in $M$ and $\lambda$ in $\Z_p^\times$.  

\begin{theo} \label{fullpre}
Let $G$ be a closed subgroup of $\SP_{2g}(\Z_p)$ containing transvections.  If $\ov{G}$ is irreducible and generated by transvections, then $G$ is saturated in $\SP_{2g}(\Z_p)$.
\end{theo}

This assertion is well-known when $\ov{G} = \SP_{2g}(\F_p)$, so we need only consider $p = 2$ and proper subgroups, thanks to a classical result of McLaughlin. 

\begin{prop}[\cite{McL}]   \label{Mc}
For $g \ge 2$, let $H$ be an irreducible proper subgroup of $\SP_{2g}(\F_p)$ generated by transvections.  Then $p = 2$ and $H$ is one of the following:
\begin{enumerate}[{\rm i)}]
\item  the symmetric group $\cS_m$ with $m =2g+1$ or $2g+2$,  \vspace{2 pt}

\item  the orthogonal group $O^+_{2g}(\F_2)$ with $g \ge 3$, or $O^-_{2g}(\F_2)$.
\end{enumerate}  
\end{prop}

Orthogonal groups and theta characteristics are reviewed in \cite{Dye,GrHa, BK2}.   Set 
$
\fraksp_{2g}(\F_p) = \{{A} \in {\rm Mat}_{2g} (\F_p) \, \vert \, {A}^t \, {J} + {J} \,{A}  = 0 \},
$ 
where $J$ is the {\em Gram matrix} of a basis $e_1, \dots, e_{2g}$ for $M.$ 
Then $\dim_{\F_p} \fraksp_{2g}(\F_p)  = 2g^2 + g$, since 
\begin{equation} \label{fraksp}
\fraksp_{2g}(\F_p) = \left\{ A  = \left[ \begin{smallmatrix} a&b\\c&d\end{smallmatrix} \right] \, \vert \hspace{2 pt} \, b = b^t, c = c^t, d=-a^t  \right\} \, \text{ when } \, J= \left[ \begin{smallmatrix} \hfill 0_g&I_g\\-I_g&0_g\end{smallmatrix} \right].
\end{equation}

To prove the Theorem, one verifies that sending an element $1+p A$ of $G^{(1)}$ to $A \bmod{p}$ induces an isomorphism $\cL\!: \, G^{(1)}/G^{(2)} \to \fraksp_{2g}(\F_p)$.  It follows that $1+p^n A \mapsto A \bmod{p}$ gives an isomorphism $G^{(n)}/G^{(n+1)} \to \fraksp_{2g}(\F_p)$ for all $n \ge 1.$   Then one shows by induction that $G^{(1)} \to \SP^{(1)}_{2g}(\Z/p^n \Z)$ is surjective and passes to the limit.  For the groups in the Theorem, the transvections $s$ in $\ov{G}$ form one conjugacy class.  Since $G$ contains a transvection, $s$ lifts to a transvection $\sigma$ in $G$, say $\sigma(x) = x - \lambda \, [y,x] \, y$.  Furthermore, $\cL(\sigma^p) \in \fraksp_{2g}(\F_p)$ is the matrix representation of the endomorphism of $\ov{M} = M/pM$ given by
\begin{equation} \label{LofTrans}
f_{\ov{y}}\!: \, \ov{x} \mapsto (s-1)(\ov{x}) =  - \lambda \, [\ov{y},\ov{x}] \, \ov{y},
\end{equation}
where $x \mapsto \ov{x}$ denotes the projection map $M \to \ov{M}$ and the pairing is induced on $\ov{M}$.  Hence it suffices to show that the $f_{\ov{y}}$'s generate $\fraksp_{2g}(\F_p)$.  This is done in Lemmas \ref{fspan} and \ref{orthog}.

\subsection{The case $\ov{G} = \cS_m$} Let $\cS_m$ act by permuting the coordinates of
$$
W=\{(a_1,..,a_m) \in \, \F_2^m\ |\  a_1 + \dots + a_m =0\}
$$
and let $V=W/\{ (a,\dots,a) \, | \, a \in \F_2 \}$ or $V=W$ if $m$ is even or odd, respectively.   Then $\ov{G}$ is symplectic for the pairing on $V$ induced by $[(a_i),(b_i)] =\sum_1^m a_i  b_i$. 
 
Let $v_{ij}$ in $V$ be represented by the vector in $W$ with non-zero entries only in coordinates $i$ and $j$.   For $x$ in $V$, the transvection $s_{ij}(x) = x + [v_{ij},x] \, v_{ij}$ corresponds to $\pi = (ij)$ in $\cS_m$.   Let $f_{ij}(x) = (s_{ij}-1)(x) = [v_{ij},x] \, v_{ij}$ be the endomorphism  of $V$ as in \eqref{LofTrans}.   

\begin{lem}  \label{fspan}
The set $\{f_{ij} \, \vert \, 1 \le i < j \le 2g+1 \}$ spans $\fraksp_{2g}(\F_2)$.
\end{lem}

\begin{proof}
The vectors $b_i = v_{i,2g+1}$ for $1 \le i \le 2g$ form a basis for $V$.  By definition, $v_{ij} = b_i+b_j$ for $1 \le i < j \le 2g$.  To prove that the $2g^2+g$ elements in the Lemma span, we show that they are linearly independent.  If not, there are constants $\alpha_{ij}, \beta_i \in \F_2$ satisfying 
$$
\sum_{1 \le i < j \le 2g} \alpha_{ij} [x,v_{ij}] \, v_{ij} \hspace{3 pt} +  \hspace{4 pt} \sum_{i=1}^{2g} \,  \beta_i [x,b_i] \, b_i = 0 \hspace{15 pt} \text{for all } x \in V.  \vspace{4 pt}
$$ 
Fix $k$ in $\{1, \dots, 2g\}$ and evaluate at $x = b_k$, using
$$
[b_k, v_{ij}] = \begin{cases} 1 &\text{if } i = k \text{ or } j = k, \\ 0 &\text{otherwise}\end{cases}
\hspace{20 pt}  \text{and} \hspace{20 pt} 
[b_k, b_i] = \begin{cases} 1 &\text{if } i \ne k, \\ 0 &\text{if } i = k, \end{cases} 
$$
to obtain
\begin{equation} \label{eq0}
\sum_{ j > k} \alpha_{kj} (b_k + b_j) + \sum_{ i < k} \alpha_{ik} (b_i + b_k)+ \sum_{i \ne k} \beta_i b_i = 0.
\end{equation}  
Match coefficients of $b_i$ with $i < k$ to find that 
$\alpha_{ik} = \beta_i $ for all $ i < k
$ 
and those  of $b_j$ with $j > k$ to find that 
$
\alpha_{kj} = \beta_j $ for all $ j >  k.
$
Hence $\alpha_{ij} = \beta_i = \beta_j = \gamma$ is constant.  From the coefficients of $b_k$ in (\ref{eq0}) we have $(2g-1) \gamma b_k = 0$ and so $\gamma = 0.$
\end{proof}

\subsection{The case $\ov{G} \simeq O^{\pm}_{2g}(\F_2)$}  
Let $V$ be a symplectic space of dimension $2g$ over $\F_2$ with basis $\{b_i\}$ and Gram matrix $J$ in \eqref{fraksp}. Consider the theta characteristic
\begin{equation} \label{theta}
 \theta_\epsilon(x) = Q_{\epsilon}(x_g,x_{2g}) + \sum_{i=1}^{g-1} x_j x_{g+j} \hspace{2 pt} \text{ for } \hspace{2 pt} x = (x_1, \dots, x_{2g}) \text{ in } V,
\end{equation}  
with  $Q_{+}(x_g,x_{2g}) = x_g x_{2g}$ in the even case and $Q_{-}(x_g,x_{2g}) = x_g^2 + x_gx_{2g} + x_{2g}^2$ in  the odd case.  We have $\theta_\epsilon(x+y) = \theta_\epsilon(x) + \theta_\epsilon(y) + [x,y]$, i.e.\! $\theta_\epsilon$ belongs to the pairing $[ \hspace{3 pt}, \hspace{3 pt}]$ on $V$.   The transvection $s\!:\, x \mapsto x + [v,x] \, v$ in $\SP_{2g}(\F_2)$ acts on theta characteristics by  
$
s(\theta)(x) = \theta(x) + (1 + \theta(v))[v,x]^2.
$ 
Hence $s$ is in the stabilizer $O_{2g}^\epsilon(\F_2)$ of $\theta_\epsilon$ in $\SP_{2g}(\F_2)$ exactly when $\theta_\epsilon(v) = 1$.  

Let   $f_{v}\!:x\mapsto [v,x] \,v$ be  the endomorphism  of $V$ as in \eqref{LofTrans}. For $\epsilon = \pm$, let
$$
\cF_g^\epsilon(V) = \{ f_v \in \fraksp_{2g}(\F_2) \, \vert \, v \in V \text{ and } \theta_\epsilon(v) = 1 \}
$$
and write $\lr{\cF_g^\epsilon(V) }$ for its span in $\fraksp_{2g}(\F_2)$.

\begin{Rem} \label{smallg} $\lr{\cF_1^-(V)} = \fraksp_2(\F_2) = \{ A \in \Mat_2(\F_2) \, \vert \, {\rm trace}(A) = 0\}$, with 
$$
f_{b_1} =  \left[ \begin{smallmatrix} 0 & 1 \\ 0 & 0 \end{smallmatrix} \right], \quad  f_{b_2} =  \left[ \begin{smallmatrix} 1 & 0 \\ 0 & 0 \end{smallmatrix} \right] \quad \text{and} \quad f_{b_1+b_2} =  \left[ \begin{smallmatrix} 1 & 1 \\1 & 1 \end{smallmatrix} \right].
$$
Also, $\dim \lr{\cF_1^+(V)} = 1$,  $\dim \lr{\cF_2^+(V)} = 6$ and $\lr{\cF_3^+(V)} = \fraksp_6(\F_2)$ is 21-dimensional.  There are symplectic isomorphisms $O_6^+(\F_2) \simeq \cS_8$ and $O_4^-(\F_2) \simeq \cS_5$.  
\end{Rem}

\begin{lem}  \label{orthog}
We have $\lr{\cF_g^\epsilon(V)} = \fraksp_{2g}(\F_2)$ \, for \,  $\begin{cases} g \ge 1 &\text{if }\epsilon =-, \\ g \ge 3 &\text{if } \epsilon = +. \end{cases}$
\end{lem}

\begin{proof}
For $g \ge 2$, the subspace $V_1 = \spn\{b_j \, \vert \, j \ne 1,g+1\}$ of $V$ is isomorphic to the symplectic space of dimension $2g-2$ whose theta characteristic  $\theta'_\epsilon$ and pairing are obtained by restriction from $\theta_\epsilon$ and $[ \hspace{3 pt}, \hspace{3 pt}]$.  Given a linear map $h_1\!: \, V_1 \to V_1$, let $h = \delta_1(h_1)$ denote its extension to $V$ satisfying $h(b_1) =  h(b_{g+1}) = 0$.  Then $\delta_1\!:  \, \fraksp_{2g-2}(\F_2) \hookrightarrow \fraksp_{2g}(\F_2)$.  In particular, $\delta_1(x' \mapsto [v', x'] \, v')$ is the matrix of the linear map $f_{v'}$ on $V$ given by $f_{v'}(x) =  [v', x] \, v'$.  Hence $Y_1 = \delta_1(\cF_g^\epsilon(V_1))$ is contained in $\cF_g^\epsilon(V)$.   

Remark \ref{smallg} treats the small values of $g$, so we assume $g \ge 2$ for $\cF_g^-(V)$ and $g \ge 4$ for $\cF_g^+(V)$.  Let $V_2 = \spn\{e_j \, \vert \, j \ne 2,g+2\}$ and extend $h_2\!:  V_2 \to V_2$ to a linear map $h = \delta_2(h_2)$ on $V$ by setting $h(b_2) =  h(b_{g+2}) = 0$.  Then $\delta_2\!:  \, \fraksp_{2g-2}(\F_2) \hookrightarrow \fraksp_{2g}(\F_2)$ and $Y_2 = \delta_2(\cF_g^\epsilon(V_2))$ also is contained in $\cF_g^\epsilon(V)$.  By the induction hypothesis, 
$$
\dim Y_i = \dim \fraksp_{2g-2}(\F_2) = 2(g-1)^2+(g-1) = 2g^2-3g+1 \quad \text{for } i = 1,2.
$$  
We have $\dim Y_1 \cap Y_2  \le \dim \fraksp_{2g-4}(\F_2) = 2(g-2)^2+g-2 = 2g^2-7g+6$.   But $\dim \fraksp_{2g}(\F_2) = 2g^2+g$, so the codimension of $Y_1+Y_2$ in $\fraksp_{2g}(\F_2)$ is at most 4. 

Next we fill in a 4-dimensional subspace of $\fraksp_{2g}(\F_2)$  independent of $Y_1+Y_2$.  For matrices in $Y_k$, all entries in the row or column numbered $k$ or $g+k$ are 0.  If $A$ is in $Y_1+Y_2$, then $A_{ij} = 0$ for the eight pairs $(i,j)$ with 
$$
\{i,j\} = \{1,2\}, \{1,g+2\}, \{g+1,2\} \text{ or } \{g+1,g+2\}.
$$
By the symplectic condition \eqref{fraksp} on $A$, we have 
\begin{equation*}
A_{1,2} = A_{g+2,g+1}, \hspace{8 pt} A_{1,g+2} = A_{2,g+1}, \hspace{8 pt}   A_{g+1,2} =  A_{g+2,1}, \hspace{8 pt}   A_{g+1,g+2} =  A_{2,1}. 
\end{equation*}
Define $j^* = j+g$ if $j < g$ and $j^* = j-g$ otherwise.  Fix $(i,j)$ in 
$$
S = \{ (1,2), \, (1,g+2), \, (g+1,2), \, (g+1,g+2) \}
$$
 and let $v = v_{ij} = b_i +b_{j^*} + b_g+b_{2g}$.  Since $\theta_\epsilon(v) = 1$, the linear map $f_v$ is in $\cF_g^\epsilon(V)$.  Also, $A_{i,j} = 1$ and $A_{i',j'} = 0$ for all other pairs $(i',j')$ in $S$.  Hence the $f_v$'s generate the desired 4-dimensional space as $(i,j)$ ranges over the 4 pairs in $S$.  
\end{proof}

In the following Proposition, a group  is said to be McL if it is isomorphic to $\SP_{2g}(\F_2)$ or one of the groups in  Proposition \ref{Mc}. 
\begin{prop}  \label{BigGal}
Let $A_{/\Q}$ be a $g$-dimensional abelian variety of odd conductor $N$, with $q\Vert N$ for a prime $q$ ramifying in $F = \Q(A[2])$ and let $F_\infty = \Q(A[2^\infty])$.  If $\ov{G} = \Gal(F/\Q)$ is {\rm McL}, then $G = \Gal(F_\infty/\Q)$ is saturated in $\GS_{2g}(\Z_2)$.  Also, $\ov{G}$ is {\rm McL} if $A[2]$ is irreducible and one of the following holds:
\begin{enumerate}[{\rm i)}]
\item the conductor of $A[2]$ is square-free and $\sqrt{-1}$ is not in $F$, or 
\item $F$ contains no proper extension of $\Q$ unramified at $q$.
\end{enumerate}
\end{prop}

\begin{proof} 
Since $A[2]$ is irreducible in all cases, a minimal polarization of $A$ has odd degree.  Hence the Weil pairing induces a perfect pairing on the  Tate module $\T_2(A)$ and $G$ is a closed subgroup of $\GS_{2g}(\Z_2)$.  The symplectic similitude $\nu\!: G \to \Z_2^\times$ giving the action of $G$ on $\Mu_{2^\infty}$ is surjective and so $F_\infty$ contains $\Q(\Mu_{2^\infty})$.  By Grothendieck's monodromy theorem, inertia at $v \vert q$ is generated topologically by a transvection $\sigma_v,$ since the toroidal dimension of $A$ at $q$ is 1 and $q$ ramifies in $F$.  

Assume (i).  For each prime $w$ dividing cond$(A[2])$, inertia at $w$ is generated by a transvection $s_w$ in $\ov{G}$. Thus the fixed field $k$ of the normal subgroup generated by all $s_w$ is unramified outside 2$\infty$.   As in \cite[Prop.\! 6.2]{BK2}, Fontaine's bound on ramification at 2 implies that $k$ is contained in $\Q(i)$ and thus $k = \Q$.  Hence $\ov{G}$ is generated by transvections and so $\ov{G}$ is McL.  In case (ii), the subfield of $F$ fixed by the normal closure of $\ov{\sigma}_v$ in $\ov{G}$ is unramified over $q$, so equals $\Q$ and $\ov{G}$ is McL.

If $\ov{G}$ is McL, transvections form one conjugacy class generating $\ov{G}$.  Since $\sigma_v$ fixes $F \cap \Q(\Mu_{2^\infty})$, the latter equals $\Q$ and the restriction of $\nu$ to $G^{(1)} = \Gal(F_\infty/F)$ surjects onto $\Z_2^\times$.  For $H = G \, \cap \,  \SP_{2g}(\Z_2)$, we find that $\ov{H} = \ov{G}$ so $H^{(1)}  = \SP^{(1)} _{2g}(\Z_2)$ by Theorem \ref{fullpre}.  Hence $G^{(1)} = \GS^{(1)} _{2g}(\Z_2)$.
\end{proof}

\begin{Ex} 
Let $A$ be the Jacobian of a hyperelliptic curve $y^2 + Q(x)y = P(x)$, where $Q(x)^2 + 4P(x)$ is irreducible in $\Z[x]$ of degree $m=2g+1$ or $2g+2$, with Galois group $\cS_m \subset \SP_{2g}(\F_2)$.  If the Igusa discriminant $I_{10}$ of $C$ is odd and some prime $q$ exactly divides $I_{10}$, then $G = \Gal(\Q(A[2^\infty])/\Q)$ is saturated in $\GS_{2g}(\Z_2)$.  
\end{Ex}

\section{Preliminaries on Honda systems}   \label{Prelim} \numberwithin{equation}{section} 

The  basic material on Honda systems may be found in \cite{BrCo,Con2,Fon1} and is summarized in \cite{BK3}.  We review the required notation and recall some finite Honda systems constructed in \cite{BK3}. 

Let $k$ a perfect field of prime characteristic $p$ and $\W=\W(k)$ the ring of Witt vectors over $k$.  Fix an algebraic closure $\ov{K}$ of the field of fractions $K$ of $\W$, let $\ov{\W}$ be its ring of integers and write $G_K = \Gal(\ov{K}/K)$.  
Let $\sigma\!: \, \W \to \W$  be the Frobenius automorphism characterized by $\sigma(x)\equiv x^p  \pmod{p}$ for $x$ in $\W$. The Dieudonn\'e ring $D_k=\W[\rF,\rV]$ is generated by the Frobenius operator $\rF$ and Verschiebung operator $\rV$, with $\rF\rV=\rV\rF=p$, $\rF a=\sigma(a)\rF$ and $\rV a=\sigma^{-1}(a)\rV$ for all $a$ in $\W$.   

A {\em Honda system} consists of a finitely generated free $\W$-module $\cM$, a submodule $\cL$ of $\cM$ and a Frobenius-semilinear injective endomorphism $\rF\!$ of $\cM$ such that $p\cM \subseteq \rF \cM$ and inclusion induces an isomorphism $\cL/p\cL \to \cM/\rF\cM$.   Then $\cM$ becomes a $D_k$-module with {\em Verschiebung} defined by $\rV x = \rF^{-1}(px)$ for all $x$ in $\cM$.  

Let $\cM$ be a $D_k$-module, finitely generated and free as a $\W$-module and let $\cL$ be a $\W$-submodule of $\cM$.  Then $(\cM,\cL)$ is a Honda system if and only if the following sequence is exact:
\begin{equation} \label{HondaSeq}
 0 \to \cL \xrightarrow{\rV} \cM \xrightarrow{\rF} \cM/\cL \to 0.
\end{equation}

\begin{lem}  \label{dual}
Given a Honda system $(\cM,\cL)$, let $\cM^* = \Hom_\W(\cM,\W)$ and let $\cL^*$ be the annihilator of $\cL$ in $\cM^*$.  Define $\rF$ and $\rV$ on elements $\psi$ of $\cM^*$ by
\begin{equation} \label{FVstar}
\rF(\psi)(x)  = \sigma(\psi(\rV x)) \quad \text{and}  \quad \rV(\psi)(x) = \sigma^{-1}(\psi(\rF x)) \quad \text{for all } x \in \cM.
\end{equation}
Then $(\cM^*,\cL^*)$  forms a Honda system.
\end{lem}
\begin{proof} Since $\cL \cap \rF \cM = p \cL$, the quotient $\cM/\cL$ is torsion-free, so $\cL$ is a direct summand of $\cM.$ The pairing $\lr{-,-}: \, \cM^* \times \cM \to \W$ induces perfect pairings:
$$
\cM^*/\cL^*  \times \cL \to \W \quad \text{and} \quad  \cL^* \times \cM/\cL \to \W.
$$
By dualizing \eqref{HondaSeq}, the sequence $0 \to \cL^* \xrightarrow{\rV} \cM^* \xrightarrow{\rF} \cM^*/\cL^* \to 0$ is exact.
\end{proof} 

If $\rF$ is topologically nilpotent, then $(\cM,\cL)$ is {\em connected}.    If both $F$ and $\rV$ are topologically nilpotent, then $(\cM,\cL)$ is {\em biconnected}.  

A {\em finite Honda system} is a pair $(\rM,\rL)$ consisting of a $D_k$-module $\rM$ of finite $\W$-length and a $\W$-submodule $\rL$ with $\rV\!: \, \rL\to \rM$  injective and the map $\rL/p\rL\to \rM/\rF\rM$  induced by the identity is an isomorphism.  If $(\cM,\cL)$ is a Honda sytem then $(\cM/p^n\cM,\cL/p^n\cL)$ is a finite Honda system.

Let $\widehat{\CW}_k$ denote the formal $k$-group scheme associated to the {\em Witt covector} group functor $\CW_k$, cf.\! \cite{Con2,Fon2}.  In particular, if $k'$ is a finite extension of $k$ and $K'$ is the field of fractions of $W(k')$, we have $\CW_k(k') \simeq K'/W(k')$.  If $R$ is a $k$-algebra, then $D_k = \W[\rF,\rV]$ acts on elements $\bfa = ( \dots, a_{-n}, \dots, a_{-1},a_{0})$ of $\CW_k(R)$ by $\rF \bfa = ( \dots,a_{-n}^p, \dots, a_{-1}^p,a_{0}^p)$,  $\rV \bfa = ( \dots, a_{-(n+1)}, \dots, a_{-2},a_{-1})$ and $\dot{c} \, \bfa = ( \dots, c^{p^{-n}}a_{-n}, \dots, c^{p^{-1}}a_{-1},ca_{0})$, where $\dot{c}$ is the Teichm\"{u}ller lift of $c$.  It is convenient to write $(\0,a_{-n},\dots,a_0)$ for an element of $\widehat{CW}_k(\ov{\W}/p\ov{\W})$ with $a_{-m} = 0$ for all $m > n$.

The Hasse-Witt exponential map is a homomorphism of additive groups: 
$$
\xi: \, \widehat{CW}_k(\ov{\W}/p\ov{\W}) \to \ov{K}/p\ov{\W} 
\quad \text{given by} \quad 
(\dots,a_{-n},\dots,a_{-1},a_0) \mapsto \sum \, p^{-n} \, \tilde{a}_{-n}^{p^n}
$$
independent of the choice of lifts $\tilde{a}_{-n}$ in $\ov{\W}$.  

We generally use calligraphic letters, e.g. $\cV$ for finite flat group schemes and the corresponding roman letter, e.g. $V$ for the associated Galois module.  If $(\rM,\rL)$ is the finite Honda system of $\cV$, the points of $V$ correspond to $D_k$-homomorphisms $\varphi\!: \, \rM \to  \widehat{CW}_k(\ov{\W}/p\ov{\W})$ such that $\xi(\varphi(\rL)) = 0$, with the action of $G_K$ on $V$ induced from its action on $\widehat{CW}_k(\ov{\W}/p\ov{\W})$.  

For the study of $p$-divisible groups in the next section, recall the finite Honda system $\gE_\lambda$ introduced in \cite[\S4]{BK3}  and our classification of extensions of exponent $p$ of $\gE_\lambda$ by $\gE_\lambda$.

\begin{Not}   \label{Ebasis}
Fix $\lambda$ in $k^\times$ and let $\gE_{\lambda} = (\rM,\rL)$ be the finite Honda system with a {\em standard} $k$-{\em basis} $x_1,x_2,x_3,x_4$ for $M$ such that $\rL = \spn\{x_1,x_2\}$ and Verschiebung and Frobenius are represented by the matrices: 
\begin{equation*} 
\rV=\left[\begin{smallmatrix}0&0&0&0\\ 1&0&0&0\\0&\lambda&0&0\\ 0&0&0&0\end{smallmatrix}\right], \quad \rF=\left[\begin{smallmatrix}0&0&0&0\\ 0&0&0&0\\0&0&0&1\\ 1&0&0&0\end{smallmatrix}\right].
\end{equation*}
Denote the corresponding group scheme by $\cE_\lambda$ and its Galois module by $E_\lambda$.    
\end{Not}

\begin{prop}[{\cite[Prop.~5.1.1]{BK3}}] \label{FieldForE} 
Let $\dot{\lambda}$ be the Teichm\"{u}ller lift of $\lambda$ and let $\gR_\lambda = \{ a \in \ov{\W}/p \ov{\W} \hspace{3 pt} \vert \hspace{3 pt} \lambda^{p^2} a^{p^4} \equiv  (-p)^{p+1}a \pmod{p^{p+2} \ov{\W}} \}$. For $a$ in $\gR_\lambda$, define $b=b_a$ and $c=c_a$ in $\ov{\W}/p \ov{\W}$ by $b \equiv - \frac{1}{p}\lambda^p a^{p^3} \pmod{p\ov{\W}}$ and $c \equiv \lambda a^{p^2} \pmod{p\ov{\W}}$.
\begin{enumerate}[{\rm i)}]
\item Let $x_1, \dots, x_4$ be a standard basis for the finite Honda system $\gE_\lambda$ of $\cE_\lambda$.  A $D_k$-map $\psi$ represents a point of $\cE_\lambda$ if and only if $\psi(x_1) =  (\0,c,b,a)$ for some $a$ in $\gR_\lambda$.  If so, $\psi(x_2)=(\0,c,b)$, $\psi(x_3)=(\0,\lambda^{-1} c)$ and $\psi(x_4)=(\0,a^p)$.   \vspace{2 pt}

\item $F = K(E_\lambda)$ is the splitting field of $\dot{\lambda}^{p^2}x^{{p^4}-1} - (-p)^{p+1}$ over $K.$  The maximal subfield of $F$ unramified over $K$ is $F_0 = K(\Mu_{{p^4}-1}, \xi)$, where $\xi$ is any root of $x^{p+1} - \dot{\lambda}$.  Moreover $F/F_0$ is tamely ramified of degree $t = (p^2+1)(p-1)$.  For $a \ne 0$ we have
$
\textstyle{\ord_p(a) = \frac{1}{t}, \quad \ord_p(b) = \frac{p^2-p+1}{t}, \quad \ord_p(c) = \frac{p^2}{t}.}
$
\item   $\gR_\lambda$ is an $\F_{p^4}$-vector space under the usual operations in $\ov{\W}/p \ov{\W}$ and $a \mapsto P_a$ defines an $\F_p[G_K]$-isomorphism $\gR_\lambda \xrightarrow{\sim} E_\lambda$.     
\end{enumerate}
\end{prop}

Let $\Ext^1(\gE_\lambda, \gE_\lambda)$ be the group of classes of extensions of finite Honda systems:
\begin{equation} \label{Extension}
0 \to \gE_\lambda \xrightarrow{\iota} (\rM,\rL) \xrightarrow{\pi} \gE_\lambda \to 0
\end{equation}
under Baer sum.  The subgroup $\Ext^1_{[p]}(\gE_\lambda, \gE_\lambda)$ of those classes such that $p \rM = 0$ was determined in \cite[Prop.~4.5]{BK3}, as follows.

\begin{prop} \label{Hexpp}
If $(\rM,\rL)$ represents a class in $\Ext^1_{[p]}(\gE_\lambda, \gE_\lambda)$, then there is a $k$-basis $e_1, \dots, e_8$ for $\rM$ such that $\iota(x_1) = e_1$, $\pi(e_5) = x_1$, $\rL = \spn\{ e_1,e_2,e_5,e_6\}$,   
\vspace{2 pt}
{\small 
$$
\rV=\left[\begin{array}{cccc|cccc}
              0&0&0&0&0&\lambda s_2&0&0 \\
              1&0&0&0&0&\lambda s_3&0&0 \\
              0&\lambda&0&0&0&\lambda s_4 &0&0  \\
              0&0&0&0&s_1&\lambda s_5&0&0 \\
                \hline
              0&0&0&0&0&0&0 &0   \\ 
              0&0&0&0&1&0&0&0        \\ 
              0&0&0&0&0&\lambda&0&0    \\ 
              0&0&0&0&0&0&0&0                                                                                
                                             \end{array}\right] \text{and } \,
\rF= \left[\begin{array}{cccc|cccc}0&0&0&0&0&0&0&0
                                            \\ 0&0&0&0&0&0&0&0
                                             \\ 0&0&0&1&0&-s_1^p&-s_5^p&0
                                             \\ 1&0&0&0&0&0&-s_2^p&0
                                             \\ \hline 0&0&0&0&0&0&0&0
                                             \\ 0&0&0&0&0&0&0&0
                                             \\ 0&0&0&0&0&0&0&1
                                             \\ 0&0&0&0&1&0&0&0   \end{array}\right] 
$$}

\noindent    with $s_1, s_2, s_3, s_4, s_5$ in $k$.   For $\tilde{k} = k/(\sigma^4-1)(k)$, the map $(\rM,\rL) \leadsto (s_1, \dots, s_5)$ induces an isomorphism of additive groups
$
\bfs\!: \, \Ext^1_{[p]}(\gE_\lambda, \gE_\lambda) \xrightarrow{\sim} k \oplus k  \oplus k \oplus  \tilde{k}  \oplus k.
$ 
\end{prop}

\begin{prop}[{\cite[Prop.~5.2.16]{BK3}}] \label{CondExpP}
Let $\cW$ be an extension of $\cE_\lambda$ by $\cE_\lambda$ killed by $p$.  The field of points $L = K(\cW)$ is an elementary abelian $p$-extension of $F = K(\cE_\lambda)$ whose conductor exponent satisfies $\gf(L/F) \le p^2$.  
\end{prop}

\section{Our $p$-divisible groups}   \label{pdiv} \numberwithin{equation}{section}
We classify Honda systems $(\cM,\cL)$ associated to $p$-divisible groups whose first layer is $\cE_\alpha$, as in Notation \ref{Ebasis}, with  $\alpha$ in $k^\times$.    We also determine the Honda systems of the Cartier duals of such $p$-divisible groups.

\begin{prop}  \label{ourHons} 
Let $(\cM,\cL)$ be a Honda system as above. Then there is a basis $e_1, e_2, e_3, e_4$ for $\cM$ over $\W$ and parameters $\lambda$ in $\W^{\times}$ and $s_1,s_2,s_3,s_5$ in $\W$ such that $\lambda \equiv \alpha \pmod{p\W}$,  \, $\cL = \spn \{e_1,e_2 \}$, 
$$
\rV=\left[\begin{array}{cccc}0&p\lambda s_2&0&p\\ 1&p\lambda s_3&0&0\\0&\lambda&0&0\\ ps_1&p\lambda s_5&p&0\end{array}\right] \quad \text { and } \quad \rF=\sigma\left[\begin{array}{cccc}0&p&-p^2 s_3& 0 \\ 0&0&p/\lambda&0\\0&-p s_1 &p^2 s_1 s_3-ps_5 & 1\\ 1&0&-p s_2 & 0\end{array}\right]\!.
$$
\end{prop}

\begin{proof}
Choose a lift $u$ in $\W^\times$ of $\alpha^{-1}$ and lifts $e_{i,1}$ in $\cM$ of a standard basis for $\gE_{\alpha}$, such that $e_{4,1} = \rF e_{1,1}$ and $e_{3,1} = \rF e_{4,1}$.  We prove by induction that there is a basis $e_{1,n}, \dots, e_{4,n}$ for $\cM$ satisfying $e_{4,n} = \rF e_{1,n}$, \, $e_{3,n} = \rF e_{4,n}$, 
\begin{equation} \label{indhyp}
e_{2,n} - \rV e_{1,n}  \in  p^n \, \spn\{e_{3,n}\} + p \, \spn\{e_{4,n}\} \hspace{2 pt} \text{ and } \hspace{2 pt}e_{3,n} \in  u \rV e_{2,n} + p\cM . 
\end{equation}
Substituting the last relation into the first relation in \eqref{indhyp}, we get 
\begin{eqnarray*}
e_{2,n} - \rV e_{1,n} & \in &  p^n (a u \, \rV e_{2,n} + p\cM ) + p \, \spn\{e_{4,n}\}  \\
                   & \in &  p^n a u \, \rV e_{2,n} + p^{n+1} \cM  + p \, \spn\{e_{4,n}\} .
\end{eqnarray*}
Let $e_{1,n+1} = e_{1,n} - \sigma(a u) p^n e_{2,n}$.    Then
\begin{equation} \label{ind1}
 e_{2,n} - \rV e_{1,n+1}  \in  p^{n+1} \cM  + p \, \spn\{e_{4,n}\}. 
\end{equation}
Set $e_{4,n+1} = \rF e_{1,n+1}$.  We have  $\rF e_{2,n}\in p\cM $ by  Notation \ref{Ebasis} and $\rF e_{1,n} = e_{4,n}$, so
$$
e_{4,n+1} = \rF e_{1,n+1}  \text{ is in } \rF e_{1,n} + p^{n} \, \spn\{\rF e_{2,n}\} \subseteq e_{4,n} + p^{n+1} \cM .
$$
Hence (\ref{ind1}) is equivalent to 
\begin{equation} \label{ind2}
e_{2,n} - \rV e_{1,n+1}  \in   p^{n+1} \cM  + p \, \spn\{e_{4,n+1}\}. 
\end{equation}
Define $e_{3,n+1} = \rF e_{4,n+1}$ and use $\rF e_{4,n} = e_{3,n}$ to obtain
$$
e_{3,n+1} = \rF e_{4,n+1}  \in  \rF e_{4,n} + p^{n+1} \rF \cM  \subseteq e_{3,n} + p^{n+1} \cM .
$$
Clearly $\cM  = \spn\{e_{1,n}, \, e_{2,n}, \, e_{3,n+1}, \,e_{4,n+1}\}$.  By (\ref{ind2}), there are scalars $b_1, b_2, b_3, b_4$ such that 
$$
e_{2,n} - \rV e_{1,n+1} = p^{n+1} \sigma(b_1) e_{1,n} + p^{n+1} \sigma(b_2) e_{2,n} + p^{n+1} b_3 e_{3,n+1} + p b_4 e_{4,n+1}. 
$$
Setting $e_{2,n+1} = e_{2,n} - p^{n+1} \sigma(b_1) e_{1,n} - p^{n+1} \sigma(b_2) e_{2,n}$ gives the induction step for the first relation in (\ref{indhyp}).  For the second part,  
$$
e_{3,n+1} \text{ is in }e_{3,n} + p^{n+1} \cM  \subseteq u \rV e_{2,n} + p\cM  \subseteq u \, \rV e_{2,n+1} + p\cM .
$$
Then $e_{1,n}$, $e_{2,n}$, $e_{3,n}$, $e_{4,n}$ converge to a basis $e_1$, $e_2$, $e_3$, $e_4$ for $\cM$.  By the first part of \eqref{indhyp}, $\rV e_1 = e_2 + ps_1e_4$ for some $s_1$ in $\W$.  By the second part,  $\rV  e_2$ is in $\lambda e_3 + p \, \spn\{e_1,e_2,e_4\}$ for some $\lambda$ in $\W^\times$ lifting $\alpha$.  Also, $\rV  e_3 = \rV \rF e_4 = p e_4$ and $\rV  e_4 = \rV \rF e_1 = p e_1$.  This verifies the matrix for $\rV $ and that for $\rF  = p\rV ^{-1}$ follows by  semi-linearity.
\end{proof} 

\begin{Def} \label{HParam}
 A basis $\cB = \{e_1,e_2,e_3,e_4\}$ as in the Proposition is  a {\em standard basis} for $(\cM,\cL)$.   A standard basis for a finite Honda system $(\rM,\rL) = (\cM/p^n \cM,\cL/p^n \cL)$ is the reduction of a standard basis for $(\cM,\cL)$ viewed over $\W/p^n \W$.  Denote the associated parameters by $\bfs_\cB = [\lambda;s_1,s_2,s_3,s_5]$. 
\end{Def}

\begin{cor} \label{AutP2}
Another basis $\cB' = \{e'_i\}$ for $(\cM,\cL)$ is  standard  if and only if there is an $a$ in $\W^\times$ such that $e_1' = \sigma^2(a) e_1$,  $e_2' =\sigma(a)e_2$, $e_3' = \sigma^4(a) e_3$ and $e_4'=\sigma^3(a)e_4.$   Then $\bfs_{\cB'}$ is given by
$$
\lambda' = \frac{a}{\sigma^4(a)} \lambda, \hspace{4 pt}  \, s_1' = \frac{\sigma(a)}{\sigma^3(a)}s_1,  \hspace{4 pt}s_2' = \frac{\sigma^4(a)}{\sigma^2(a)} s_2, \hspace{4 pt} s_3' = \frac{\sigma^4(a)}{\sigma(a)} s_3 \, \text{ and } \, s_5' = \frac{\sigma^4(a)}{\sigma^3(a)}s_5.
$$
\end{cor}

\proof
Since $\cL = \spn\{e_1,e_2\} = \spn\{e'_1,e'_2\}$, we can find $a,b,c,d$ in $\W$ such that
$$
e_1'= \sigma^2(a) e_1+ \sigma(b) e_2 \quad \text{and} \quad e_2' = c e_1+ \sigma(d) e_2.
$$
Then $e_4' = \rF e_1' = \sigma^3(a) e_4+ \sigma^2(b)(p e_1-p \sigma(s_1) e_3)$ and  
\begin{eqnarray*}
\rV e_1' &=& \sigma(a)\rV (e_1)+ b \rV (e_2)  \\
          &=& \sigma(a)(e_2 +ps_1e_4) + b \, (p\lambda s_2 e_1+p \lambda s_3 e_2+ \lambda e_3 + p \lambda s_5 e_4). 
\end{eqnarray*}
From the  matrix for $\rV $ on the new basis, we have:
$$
\rV e_1' = e_2' + p s_1' e_4' = c e_1 + \sigma(d) e_2 + p s_1' (\sigma^3(a)e_4 + \sigma^2(b) (p e_1 - p \sigma( s_1) e_3)).
$$
Comparing coefficients of $e_3$ in $\rV e_1'$ gives $\lambda b = -p^2 s_1' \sigma(s_1) \sigma^2(b)$, so $b = 0$ or else
$\ord_p(b) \ge 2 +  \ord_p(b)$.   By comparing coefficients of $e_1$, we find that $c = 0$.  Hence the coefficients of $e_2$ and $e_4$ give $d = a$ and $s_1' \sigma^3(a) = \sigma(a) s_1$.  We now have:  
$$
e_1' = \sigma^2(a) e_1, \quad e_2' =\sigma(a)e_2, \quad e_4'=\sigma^3(a)e_4 \quad \text{and} \quad e_3' = \rF e_4' = \sigma^4(a) e_3.
$$
Compare coefficients in $\rV e_2' = a \rV (e_2) = a \lambda (p s_2e_1+p s_3e_2+ e_3+p s_5e_4)$ and 
\begin{eqnarray*}
\rV e_2' &=& \lambda' (ps_2'e_1'+p s_3'e_2'+ e_3'+p s_5e_4' )\\
         &=& \lambda' (ps_2' \sigma^2(a) e_1 + p s_3' \sigma(a)e_2 + \sigma^4(a) e_3 + p s_5' \sigma^3(a) e_4). \quad\qed
\end{eqnarray*}

\begin{cor}  \label{DualParams}
Let $(\cM^*, \cL^*)$ be the dual of $(\cM,\cL)$ as in {\rm Lemma \ref{dual}}.  There is a standard basis $\tilde{\cB} = \{\xi_i\}$ for $\cM^*$ with $\bfs_{\tilde{\cB}} = [\lambda'; s'_1, s'_2, s'_3, s'_5]$ related to $\bfs_\cB$ by:
$$
\begin{array}{l l l l}
\lambda' = \frac{1}{\sigma^2(\lambda)}, & s'_1= -\frac{1}{\sigma(\lambda)} s_1, &  s'_2 = -\sigma^{2}(\lambda) s_2,  &s'_3= \sigma^2(\lambda) \sigma^{-1}(ps_1s_3-s_5),  \vspace{4 pt} \\ 
\multicolumn{4}{c}{s'_5 = - \sigma^2(\lambda) \sigma(s_3) -\displaystyle{\frac{ps_1\sigma^2(\lambda)}{\sigma(\lambda)}}\sigma^{-1}(ps_1s_3-s_5).}
\end{array}
$$
\end{cor}

\proof
By \eqref{FVstar} and the Proposition, the matrices for Verschiebung and Frobenius on $\cM^*$ in terms of its dual basis $e_1^*, \dots, e_4^*$ are given by
$$
\rV =\left[\begin{array}{cccc}0&0&0&1 \\ p&0&-ps_1&0 \\
                                 -p^2s_3&p/\lambda&p^2s_1s_3-ps_5&-ps_2 \\ 0&0&1&0   \end{array}\right] \hspace{-3 pt}, \hspace{6 pt}
\rF = \sigma \left[\begin{array}{cccccccc}0&1&0&ps_1  \\ 
             p\lambda s_2 &p\lambda s_3&\lambda& p\lambda s_5 \\
             0&0&0&p \\ p&0&0&0  \end{array}\right]. 
$$
Since $\cL^*$ is the annihilator of $\cL $, we have $\cL ^* = \spn\{e^*_3,e^*_4\}$.  For the standard basis, $\xi_2$ is in $\cL ^*$ while $\rF \xi_2$ is in $p\cM ^*$, so $\xi_2 = xe^*_4+pwe^*_3$ with  $x$ in $\W^\times$ and  $w$ in $\W$.  But $\{\xi_1, \xi_2\}$  is a basis for $\cL ^*$,  so $\xi_1 = ye^*_3 + ze^*_4$, with $y$  in $\W^\times$.  By scaling, assume that $y = 1$ and $\xi_1 = e^*_3 + ze^*_4$.  Apply $\rV$ and let $t_5 = p s_1 s_3 -s_5$, to find that
$$
 \rV (\xi_1) = \sigma^{-1}(z)e^*_1- ps_1e^*_2 + p[t_5-\sigma^{-1}(z)s_2]e^*_3 + e^*_4.
$$ 
Proposition \ref{ourHons} gives $\rF$ and $\rV$ on the standard basis $\tilde{\cB}$ for $\cM^*$.  Thus:
\begin{eqnarray*}
\xi_4\!\!\!\! &=&\!\!\!\! \rF(\xi_1) = p \sigma(s_1z)e^*_1+ \sigma(\lambda +p\lambda s_5z) e^*_2 + p \sigma(z)e^*_3,  \vspace{4 pt}\\ 
\rV (\xi_1)\!\!\!\! &=&\!\!\!\! \xi_2+ps'_1\xi_4 = xe^*_4+pwe^*_3+ps'_1(p \sigma(s_1z)e^*_1+ \sigma(\lambda +p\lambda s_5z) e^*_2 + p \sigma(z)e^*_3). 
\end{eqnarray*}
Equating the coefficient of $e_1^*$ gives $\sigma^{-1}(z)= p s'_1\sigma(s_1 z),$ whose valuation implies that $z= 0$.  Comparing the other coefficients gives $x=1,$ $w = t_5$ and $\sigma(\lambda) s'_1= - s_1$, so
$\xi_1=e^*_3$,\, $\xi_2 = e^*_4 + pt_5e^*_3$,\, $\xi_4 = \sigma(\lambda) e^*_2$, $\xi_3= \rF \xi_4 = \sigma^2(\lambda)(e^*_1+p\sigma(\lambda s_3)e^*_2)$.  The remaining formulas relating $\bfs_{\tilde{\cB}}$ and $\bfs_\cB$ result  from a comparison of 
\begin{eqnarray*}
\rV (\xi_2) \hspace{-4 pt} &=& \hspace{-4 pt}\, \rV (e^*_4)+p \sigma^{-1}(t_5)\rV (e^*_3) \\
\hspace{-4 pt} &=& \hspace{-4 pt} e^*_1-ps_2e^*_3+p\sigma^{-1}(t_5)(-ps_1e_2^*+ pt_5e_3^* + e_4^*)\\
\hspace{-4 pt}  &=& \hspace{-4 pt}  e^*_1- p^2s_1\sigma^{-1}(t_5)e_2^* - p(s_2 - p\sigma^{-1}(t_5)t_5) e_3^*+ p\sigma^{-1}(t_5)e_4^*
\end{eqnarray*}
with 
$$
\begin{array}{l}
 \rV (\xi_2) =\lambda' \, (ps'_2\xi_1+ps'_3\xi_2+\xi_3+ps'_5\xi_4)\\
\hspace{27 pt} = \lambda' \left( ps'_2e^*_3+ps'_3 (e^*_4 + pt_5e^*_3) +\sigma^2(\lambda)(e^*_1+p\sigma(\lambda s_3)e^*_2)+ps'_5\sigma(\lambda)e^*_2 \right).   \\
\hspace{27 pt} = \lambda' \left( \sigma^2(\lambda)e^*_1+p (\sigma^2(\lambda) \sigma(\lambda s_3)+\sigma(\lambda) s'_5)e^*_2 + p(s'_2+p t_5 s'_3)e^*_3+ps'_3 e^*_4  \right).
\end{array}
$$

\begin{cor} \label{SelfDual}
Let $\cB$ be a standard basis for $\cM$ and $\bfs_\cB = [\lambda;s_1,s_2,s_3,s_5]$.  Then $(\cM^*,\cL^*)$ and $(\cM,\cL)$ are isomorphic if and only if there are $a,b$ in $\W^\times$ such that 
$$
\lambda = -\frac{\sigma^2(a)}{a} b, \hspace{30 pt} s_5 = \frac{\sigma^3(a)}{\sigma^2(a)} \sigma(bs_3) + ps_1s_3
$$ 
and one of the following holds: \, {\rm i)} $s_1 \ne 0$ {\rm (}or $s_2 \ne 0${\rm)} and $b = 1$;\quad  {\rm ii)} $s_1 = s_2 = 0$ and $b = \pm 1$; \quad {\rm iii)} all $s_j = 0$ and $b\,  \sigma^2(b) = 1$.
\end{cor}

\begin{proof}
Assume that $(\cM,\cL)$ and $(\cM^*,\cL^*)$ are isomorphic and use the previous Corollaries for the relationship between $\bfs_{\tilde{\cB}}$ and $\bfs_{\cB}$.  In particular, $\lambda \sigma^2(\lambda) = \frac{\sigma^4(a)}{a}$ for some $a$ in $\W^\times$.  Define $b$ in $\W^\times$ by $\lambda = -\frac{\sigma^2(a)}{a} b$.  Then $b \, \sigma^2(b) = 1$ and the requirement on $s'_3$ implies our claimed formula for $s_5$.  In case (i), the condition on $s'_1$ (or $s'_2$) forces $b = 1$.  In case (ii), use $s'_3$ to find that $b = \pm 1$.  Conversely, in each of these cases, there is an isomorphism. 
\end{proof}

\begin{Ex}
Since the finite Honda system with parameters $\bfs_{\cB} = [\lambda;0,0,0,0]$ plays an important role in later conductor estimates, note that it occurs naturally for $p=2$.  Let $A = J(C)$ be the Jacobian of the curve $C\!:\,y^2 + ay = x^5 + b$ with $a$ and $b$ units in $\W.$ Then $A[2]$ is isomorphic to $\cE_\lambda$ as in Notation \ref{Ebasis}.   Without loss of generality, we may assume that a primitive fifth root of unity $\zeta_5$ is in $\W$. The automorphism $(x,y) \mapsto (\zeta_5 x, y)$, induces a complex multiplication on $A$ and so an isomorphism $\varphi$ of the associated Honda system $(\cM,\cL)$.  If $\cB$  is a standard basis, so is $\cB' = \varphi(\cB)$ and $\bfs_{\cB} = [\lambda;s_1,s_2,s_3,s_5]$ is related to $\bfs_{\cB'} = [\lambda';s'_1,s'_2,s'_3,s'_5]$ as in Corollary \ref{AutP2}.  Moreover, $\varphi(e_1) = \sigma^2(\zeta) e_1$, where $\zeta$ is a primitive fifth root of unity, since $\varphi$ has order 5 on $\cM$.   Because $\varphi$ commutes with $\rV$, each basis leads to the same matrix for $\rV$, as  in Proposition \ref{ourHons} and thus $\bfs_{\cB'} = \bfs_\cB$.  But $s'_i$ is a multiple of $s_i$ by a primitive fifth root of unity, so  each $s_i = 0$.  
\end{Ex}

\section{Exponent $p^2$}  \label{expp2} \numberwithin{equation}{subsection}

Let $(\cM,\cL)$ be the Honda system of a $p$-divisible group as in Proposition \ref{ourHons}.  Throughout this section, $\gM = (\rM,\rL)$ denotes the finite Honda system of exponent $p^2$ such that $(\rM,\rL) \simeq (\cM /p^2 \cM , \cL/p^2 \cL),$ with  standard basis $\cB = \{e_1,e_2,e_3,e_4\}$ over $\W/p^2\W$ and  parameters  $\bfs_\cB(\gM) = [\lambda;s_1,s_2,s_3,s_5]$, where $\lambda$ is in $(\W/p^2\W)^\times$ and the $s_i$ are in $k = \W/p\W$. Write $\cV$  for the  group scheme associated to $\gM.$

\subsection{Baer sums for exponent $p^2$}  \label{sec:BaerSum}
Let $x_1, \dots, x_4$ be a standard basis for $\gE_\lambda$, cf.\! Notation \ref{Ebasis}.  Then $\gM$ is an extension of $\gE_\lambda \simeq (\rM /p\rM , \rL/p \rL)$ by $\gE_\lambda \simeq (p\rM,p\rL)$:
\begin{equation} \label{ExtP2}
0 \to \gE_\lambda \xrightarrow{\iota} (\rM,\rL) \xrightarrow{\pi} \gE_\lambda \to 0,
\end{equation} 
with $D_k$-maps $\iota$ and $\pi$ induced by $\iota(x_1) = pe_1$ and $\pi(e_1) = x_1$.  Write $[\gM]$ for the class in $\Ext^1(\gE_\lambda,\gE_\lambda)$ of the extension \eqref{ExtP2}.   Let $\gM' = (\rM ',\rL')$ be the finite Honda system of Proposition {\rm \ref{Hexpp}} with $p \rM' = 0$ and parameters $\lambda \bmod{p}$ in $k^\times$ and $\bfs(\gM') = (s_1',s_2,s_3',0,s_5')$ in $k^5$ with respect to a standard basis $e_1', \dots, e_8'$.

\begin{prop}  \label{Baer}
The Baer sum $[\gM] + [\gM']$ in $\Ext^1(\gE_\lambda,\gE_\lambda)$ is the extension class of $\gM'' = (\rM'',\rL'')$, where $\rM''$ has exponent $p^2$ and its parameters with respect to a standard basis are $\bfs_{\cB''}(\gM'')  = [\lambda;s_1+s_1',s_2+s_2',s_3+s_3',s_5+s_5']$. 
\end{prop}

\begin{proof}
To construct $\gM''$, let 
$\Gamma = \{ \, (m,m')\in \rM  \oplus \rM '  \, \vert \, \pi(m) = \pi'(x') \, \}$ be the fiber product with coordinatewise action of the Dieudonn\'{e} ring $D_k$.  Impose the relations
$
\Delta=\{\, (\iota(x),0) - (0,\iota'(x)) \in \Gamma \, \vert \, x \in \gE_\lambda \, \}, 
$
so that $\rM '' = \Gamma/\Delta$ and $\rL ''$ is the image in $\rM ''$ of 
$
\{ (p,p')\in \rL  \oplus \rL '  \, \vert \, \pi(p)=\pi'(p') \}.
$
Then $\gM'' = (\rM'',\rL'')$ is an extension:
$$
0 \to \gE_\lambda \xrightarrow{\iota''} (\rM '',\rL '') \xrightarrow{\pi''} \gE_\lambda \to 0,
$$ 
with $\iota''$ induced by $x \leadsto (x,0)$ and $\pi''$ by $(m,m') \leadsto \pi(m)$.    The elements 
$$
\gamma_1 = (e_1,e_5'), \hspace{5 pt} \gamma_2= (e_2,e_6'), \hspace{5 pt} \gamma_3= (e_3, e_7'), \hspace{5 pt} \gamma_4= (e_4,  e_8')
$$ 
of $\rM  \oplus \rM '$ satisfy the fiber product condition for membership in $\Gamma$ and we claim that their cosets modulo $\Delta$ form a standard basis for $\rM ''$ over $\W/p^2\W$.  Indeed, $\rL '' = \spn\{\gamma_1,\gamma_2\}$, $\rF \gamma_1 = \gamma_4$ and $\rF \gamma_4 = \gamma_3$.  

Write  $(m_1,m_1') \sim (m_2,m_2')$ if each pair represents the same coset in $\rM''$.  The relations on  $\Gamma$ give $(0,e_i') \sim (p e_i,0)$ for $1 \le i \le 4$, so 
$$
\rV \gamma_1 = (e_2+ps_1e_4,e_6'+s_1' e_4') = (e_2, e_6') + (p(s_1+s_1') e_4,0) \sim  \gamma_2 + p (s_1+s_1') \gamma_4.
$$
 A similar computation, using  $s_4' = 0$, shows that 
\begin{equation} \label{BaerSumParam}
\rV \gamma_2 \sim \lambda (\gamma_3+ p(s_2+ s_2') \gamma_1 + p(s_3+ s_3') \gamma_2  + p(s_5+ s_5') \gamma_4). 
\end{equation}
Complete the remaining entries of $\rV$ and $\rF$  to verify that $\gamma_1, \dots, \gamma_4$ is a standard basis for $\gM''$.  Then the parameters of $\gM''$ can be read off \eqref{BaerSumParam}.
\end{proof}

\subsection{Field of points and conductor for exponent $p^2$}   \label{P2Eq}
We first determine the field of points for the group scheme $\cV^{(0)}$ associated to the finite Honda system $\gM^{(0)}$ of exponent $p^2$ with parameters $[\lambda;0,0,0,0]$.  

\begin{lem}  \label{lem:L01}
The field  $L_0 = K(\cV^{(0)})$ is the compositum of $F = K(\cE_\lambda)$ and $\Q_p(z)$, where $z$ is any root of
\begin{equation} \label{all0} 
f(Z) = \left(Z^{p^3} + \frac{1}{p}\right)^p + (-1)^p p^{p-1}Z^{p^4} -  Z - \delta_p
\end{equation}
with $\delta_p = 1$ if $p=2$ and $0$ otherwise.  
\end{lem}

\begin{proof}
Points of $\cV^{(0)}$ correspond to $D_k$-homomorphisms $\varphi\!: \, \rM \to \widehat{CW}_k(\ov{\W}/p\ov{\W})$.  Let $e_1, e_2, e_3, e_4$ be a standard basis for $\gM^{(0)} = (\rM,\rL)$.  Then $e_1$ generates $\rM$ as a $D_k$-module and $\rV ^4$ vanishes on $\rM$ by Proposition \ref{ourHons}.  Hence $\varphi$ is determined by $\varphi(e_1) = (\0,y_3,y_2,y_1,y_0)$.   By applying $\rF$ to $e_1$ and $e_4$, we find that:
\begin{equation} \label{eq:1}
\varphi(e_4)=(\0,y_3^p,y_2^p,y_1^p,y_0^p) \quad  \text{ and } \quad \varphi(e_3) = (\0,y_3^{p^2},y_2^{p^2},y_1^{p^2},y_0^{p^2}).
\end{equation}
Since $\{f_i = pe_i \, \vert \, 1 \le i \le 4\}$ is a standard basis for $(p \rM,p \rL) \simeq \gE_\lambda$, we have $\varphi(f_1) = (\0,c,b,a)$, where $a$ is a root of $\lambda^{p^2} a^{p^4-1} = (-p)^{p+1}$ as in Proposition \ref{FieldForE}.  Then 
$$
(\0,y_3^p,y_2^p,y_1^p) \, = \, \rV \rF \varphi(e_1) \, = \, \varphi(f_1) \, = \, (\0,c,b,a) 
$$
and so
$$
y_1^p\equiv a,  \hspace{15 pt}  y_2^p\equiv b\equiv -\textstyle{\frac{1}{p}} \lambda^p a^{p^3},  \hspace{15 pt}  y_3^p\equiv c\equiv \lambda a^{p^2} \pmod{p\ov{\W}}.
$$  
Also, $\varphi(e_2) = \rV \varphi(e_1) = (\0,y_3,y_2,y_1)$.   Use Proposition \ref{FieldForE}(ii) to evaluate:
$$
\ord_p(y_1)= \textstyle{\frac{1}{p}} \ord_p(a),\hspace{10 pt} \ord_p(y_2) = \textstyle{\frac{1}{p}} \ord_p(b),\hspace{10 pt} \ord_p(y_3) = \textstyle{\frac{1}{p}} \ord_p(c).
$$
With these estimates, \eqref{eq:1} gives $\varphi(e_3) = (\0,a^p,y_0^{p^2})$.     From $\rV (e_2) = \lambda e_3$, we have:   
$$
(\0,y_3,y_2) = \rV  \varphi(e_2) = \lambda \varphi(e_3) =  \lambda(\0,a^p,y_0^{p^2}) = (\0,\lambda^{1/p}a^p,\lambda y_0^{p^2}). 
$$
Thus $y_3 \equiv  \lambda^{1/p} a^p$ and $y_2 \equiv \lambda y_0^{p^2}$ $\pmod{p\ov{\W}}$. 

Vanishing of the Hasse-Witt exponential map on $\rL$ implies that:
$$
y_0+\frac{y_1^p}{p}+\frac{y_2^{p^2}}{p^2}+\frac{y_3^{p^3}}{p^3} \equiv 0
\quad\text{ and }\quad
y_1+\frac{y_2^p}{p}+\frac{y_3^{p^2}}{p^2} \equiv 0 \pmod{p\ov{\W}}
$$
and so
$$
y_0 +(-1)^p \frac{\lambda^{p^2}}{p} \left( \frac{y_0^{p^3}}{p}+\frac{a^{p^3}}{p^2} \right)^p+ \frac{\lambda^{p^2}}{p^2} y_0^{p^4} +\delta_pa  \equiv 0 \pmod{p \ov{W}}.
$$
The substitution $y_0 = az$ relates $y_0$ to a root $z$ in $\ov{\W}$ of \eqref{all0}.   Since $V^{(0)}$ is generated by the points of $E_\lambda$ and a choice of $\varphi$, we have $L_0 = F(z)$.
\end{proof}

\begin{Not}
If $x$ is an element of a $p$-adic ring $R$, the {\em big}-0 notation $0(x)$ represents an element of the ideal $(x)$.
\end{Not}

\begin{lem} \label{lem:L02}
Let $r$ be a root of $x^{p^3} - x^{p^2}+ x^p - x + \frac{1}{p}$.  Then:
\begin{enumerate}[{\rm i)}]
\item $\fdeg{\Q_p(r)}{\Q_p} = p^3$, a prime of $\Q_p(r)$ is given by $\pi = \frac{1}{r}$ and $\ord_p(r) = -\frac{1}{p^3}$.
\item The splitting field of $f(Z)$ in {\rm \eqref{all0}} is $\Q_p(\Mu_{p^4-1},r)$.
\item  The roots of $f(Z)$ have the form $r + \alpha + 0(\pi^{(p-1)p})$ with $\alpha$ in $\{0\} \cup \Mu_{p^4-1}$.
\end{enumerate}
\end{lem}

\begin{proof}
Since $\ord_p(r) = -\frac{1}{p^3}$, the equation for $r$ is irreducible over $\Q_p$ and (i) holds.  For convenience, let $\epsilon = \pi^{p^2-p}$ in $\Q_p(r)$, so $\ord_p(\epsilon) = \frac{p-1}{p^2}$.  First we show that $r$ is an approximate root of $f$ satisfying $f(r) = 0(\epsilon)$.  

Let  $s = r^{p^3} + \frac{1}{p} = r^{p^2} - r^p + r$, so $\ord_p(s) = -\frac{1}{p}$.  By binominal expansion:
$$
s^p = (r^{p^2} - r^p + r)^p = r^{p^3} +(-1)^p r^{p^2} + r^p + 0(\epsilon) =  r - \textstyle{\frac{1}{p}} + 0(\epsilon)
$$
since the worst case middle term sastisfies $\ord_p\left(pr^{p^2(p-1)}r^p\right) = \frac{p-1}{p^2}$ and also using $\ord_2(2r^4) = \frac{1}{2} > \ord_2(\epsilon)$ when $p = 2$.  Similarly, 
$$
(-1)^p p^{p-1}r^{p^4} = (-1)^p\textstyle{\frac{1}{p}} (pr^{p^3})^p = \frac{1}{p}(1-ps)^p =  \textstyle{\frac{1}{p}} + (-1)^p p^{p-1}s^p + \epsilon_1
$$
with $\ord_p(\epsilon_1) = \ord_p(ps) = \frac{p-1}{p} > \ord_p(\epsilon)$.   By the formula for $s^p$ above, we have
$$
(-1)^p p^{p-1} r^{p^4} =  \textstyle{\frac{1}{p}} + \textstyle{(-1)^p} p^{p-1} (r - \frac{1}{p} + 0(\epsilon)) + \epsilon_1  = \textstyle{\frac{1}{p}} + (-1)^{p+1} p^{p-2} + 0(\epsilon)
$$ 
since $\ord_p(p^{p-1}r) \ge 1-\frac{1}{p^3}$.    Modulo $p$, replace $(-1)^{p+1} p^{p-2}$ with $\delta_p$.  Then
\begin{eqnarray*}
f(r) &=& s^p + (-1)^p p^{p-1} r^{p^4} - r -  \delta_p \\  
         &=&  \left(r - \textstyle{ \frac{1}{p}}\right) + \left( \textstyle{ \frac{1}{p}}+ \delta_p \right) - r - \delta_p + 0(\epsilon) = 0(\epsilon).
\end{eqnarray*}

We next show that the coefficients of $g(y) = f(y+r) - f(r)$ lie in $\Z_p[\pi]$ and that $g(y) \equiv y^{p^4} - y \pmod{\pi^{p^2}}$.  If we expand
\begin{eqnarray*}
\left((y+r)^{p^3} + \textstyle{\frac{1}{p}}\right)^p &=& y^{p^4} + \left(r^{p^3} +\textstyle{\frac{1}{p}}\right)^p + h_1(y) \quad \text{and} \\
p^{p-1}(y+r)^{p^4} &=& p^{p-1} r^{p^4} + h_2(y),  \vspace{4 pt}
\end{eqnarray*}
then $g(y) = y^{p^4} - y +  h_1(y) + (-1)^p h_2(y)$.   Let $C_j(h)$ denote the coefficient of $y^j$ in any polynomial $h(y)$ and recall that $\ord_p \bino{p^n}{j} = n - \ord_p(j)$.  The terms in $h_2(y)$ involve $y^j$ for $1 \le j \le  p^4$ and we have
\begin{equation} \label{term1}
\ord_p(C_j(h_2)) = p-1 + \ord_p\left(\bino{p^4}{j} r^{p^4-j}\right) = 3 - \ord_p(j) + \frac{j}{p^3} \ge 1.
\end{equation}
(The minimum occurs for $j = p^3$; also for $j=16$ if $p=2$.)  Thus $h_2(y) \equiv 0 \bmod{p}$.  To estimate $\ord_p(C_j(h_1))$, recall that $s = r^{p^3}+\frac{1}{p}$ and $\ord_p(s) = -\frac{1}{p}$, so 
\begin{equation}  \label{badterm}
\ord_p(\bino{p}{i} s^i) = 1 - \frac{i}{p} \ge \frac{1}{p} \quad \text{for } 1 \le i \le p-1.
\end{equation}
Similarly, $(y+r)^{p^3} = y^{p^3} + r^{p^3} + h_3(y)$, where $h_3$ involves $y^j$ for $1 \le j \le p^3-1$, and we have
$$
\ord_p(C_j(h_3)) =  \ord_p\left(\bino{p^3}{j} r^{p^3-j}\right) = 2 - \ord_p(j) + \frac{j}{p^3} \, \ge \, \frac{1}{p}.
$$
(The minimum occurs when $j = p^2$.)  It follows that 
\begin{equation} \label{h3ord}
(y^{p^3} + h_3(y))^p \equiv y^{p^4} \bmod{p}.
\end{equation}
By definition of $h_1$, we now find that 
$$
\begin{array}{l l l l l}
y^{p^4} + s^p + h_1(y) &=& ((y+r)^{p^3} + \frac{1}{p})^p &=& ( y^{p^3} + r^{p^3} +  h_3(y) + \frac{1}{p})^p \\
 &=& (y^{p^3} + h_3(y) + s)^p &=&  (y^{p^3} + h_3(y))^p + s^p + h_4(y)
 \end{array}
 $$
with $\ord_p(C_j(h_4)) \ge \frac{1}{p}$ by \eqref{badterm}.  Moreover, $h_1(y) \equiv h_4(y) \bmod{p}$ by \eqref{h3ord}.  Thus $g(y) \equiv y^{p^4} - y \pmod{\pi^{p^2}}$.  

Modulo $\pi^{(p-1)p}$, we have $f(r) \equiv 0$ and thus $f(y+r) = f(r) + g(y) \equiv y^{p^4} - y$.  If $\alpha$ is in $\{0\} \cup \Mu_{p^4-1}$, then $f(r+\alpha) \equiv 0 \pmod{\pi^{(p-1)p}}$ and $f'(r+\alpha)$ is a unit.  We conclude by Hensel's Lemma that there is a root $\theta_\alpha$ of $f$ in $\Z_p[r]$ such that $\ord_p(\theta_\alpha- (r+\alpha)) \ge \frac{p-1}{p}$ and this account for all the roots of $f$.
\end{proof}

\begin{prop}  \label{CondL0}
The field of points $L_0 = F(\cV^{(0)})$ is an elementary abelian $p$-extension of $F = K(E_\lambda)$, totally ramified of degree $p^3$, with ray class conductor exponent $\gf(L_0/F) = p^3-p^2+p$.
\end{prop}

\begin{proof}
Since $H =  \Gal(L_0/F)$ is trivial on $E_\lambda$, it is an elementary abelian $p$-group.  By the Lemmas above, $L_0 = F(z) = F(r)$ because $F$ contains $\Mu_{p^4-1}$.  But the ramification in $F/\Q_p$ is tame of degree $t = (p-1)(p^2+1)$ by Proposition \ref{FieldForE}, while $\Q_p(r)/\Q_p$ is totally ramified of degree $p^3$ by Lemma \ref{lem:L02}, so $L_0/F$ also is totally ramified of degree $p^3$.  With respect to a prime element $\pi'$ of $L_0$, we have $\ord_{\pi'}(\frac{1}{r}) = p^3 t \ord_p(\frac{1}{r}) = t$.  Moreover $h(r) - r$ is a unit for all $h \ne 1$ in $H$.  By a conductor lemma \cite[A.5]{BK3}, we find that $\gf(L_0/F) = t+1 = p^3-p^2+p$.   
\end{proof}

\begin{prop}  \label{CondBd} 
Let $\cV$ be the group scheme associated to a finite Honda system $\gM = (\cM/p^2\cM, \cL/p^2\cL)$ with $(\cM,\cL)$ as in {\em Proposition \ref{ourHons}.} Then  $L = K(\cV)$ is an elementary abelian $p$-extension of $F = K(\cE_\lambda)$ with $\gf(L/F) = p^3-p^2+p$.
\end{prop}

\begin{proof}
Denote  the parameters of $\gM$ by $\bfs=\bfs_\cB(\gM) = [\lambda;s_1,s_2,s_3,s_5]$.  Proposition \ref{Baer} shows that there is a finite Honda system $\gM'$ of exponent $p$ with parameters $\lambda \bmod{p}$ and $\bfs(\gM') = (s_1,s_2,s_3,0,s_5)$ such that the Baer sum of extension classes satisfies $[\gM] = [\gM^{(0)}] + [\gM']$.   If $L' = K(\cV')$ is the field of points of the group scheme $\cV'$ associated to $\gM'$, it follows from the fiber product construction that $L \subseteq L_0L'$.  The respective conductor exponents $\gf = \gf(L/F)$, $\gf_0 = \gf(L_0/F) = p^3-p^2+p$ and $\gf' = \gf(L//F)$ satisfy $\gf \le \max\{\gf_0,\gf'\}$, cf. \cite[Lemma A.9]{BK3}.  Proposition \ref{CondExpP} asserts that $\gf' \le p^2$ and so $\gf \le \gf_0$.  But also, $[\gM^{(0)}] = [\gM] + (-[\gM'])$, so $\gf_0 \le \gf$.
\end{proof}

\begin{Rem} \label{CompareFont}
In contrast, Fontaine's bound on higher ramification implies that $\gf(L/F) \le p^3+p+1$, cf.\! \cite[Prop.~A.11]{BK3}.  The sharper bound in Proposition \ref{CondBd} is essential for our applications.  In particular, when $p =2$, we find that $\sqrt{2}$ is not in $L$, since $\gf(F(\sqrt{2})/F) = 11$.
\end{Rem} 

\subsection{Finding the Honda parameters for $A[4]$}  \numberwithin{equation}{subsection}

While it seems difficult to compute Honda parameters in general, we have the following explicit result.
 
\begin{prop}  \label{FindParams}
Let $g(x) = a_5 x^5 + a_4 x^4 + a_3 x^3 + a_2 x^2 + a_1 x + a_0$ and let $C$ be the curve $y^2 + y = g(x)$ over $\Z_2$ with $a_5$ a unit.  If $A = J(C)$ is the Jacobian of $C$, the Honda parameters associated to $A[4]$ are given by $\lambda \equiv -1 \pmod{4}$ and
$$
\begin{array}{l l l }
s_1 \equiv a_1+ a_3 a_4 +\displaystyle{\frac{a_3^2-a_3}{2}},&s_2 \equiv a_3, &s_3 \equiv s_5 \equiv a_1+a_2+a_3+a_4 \pmod{2}.
\end{array}
$$
\end{prop} 

\begin{proof}  
We sketch the argument, invoking several lemmas to be proved below.  Consider a deformation $\tilde{A} = J(\tilde{C})$, where $\tilde{C}$ is the curve $y^2 + y = \tilde{g}(x)$ with
\begin{equation} \label{tildeg}
\tilde{g} = a_5x^5 + b_3x^3+b_2x^2+b_1x+b_0.
\end{equation}
In Lemma \ref{tilder}, we prove an effective deformation result under which $\tilde{A}[4]$ and $A[4]$ share the same Honda parameters, leaving 32 explicit curves over $\Z_2$ for which the Honda parameters must be determined.   Lemma \ref{x-T} gives the Kummer group associated to $F(\frac{1}{2}P)$ for a point $P$ of order 2 on $A$.  A comparison with the Kummer group coming from the finite Honda system in \S \ref{P2Eq}, as sketched at the end of this subsection, verifies our claim
\end{proof}

\begin{Rem} Note that the eight Honda parameters $\bfs_\cB$  for $A[4]$ over $\Z_2$ are  consistent with the duality in Corollary \ref{SelfDual}.
\end{Rem}

For the lemmas needed above, let $K$ be a field not of characteristic 2. Recall the $x-T$ map (\cite{Sch}, \cite[Ch.\!\! 6]{CF}, \cite[\S 5]{PS}), which gives an explicit interpretation of the Kummer-theoretic boundary map
$
A(K)/2A(K) \hookrightarrow H^1(G_K,A[2]) 
$
arising from Galois cohomology of  
$$
0 \to A[2] \to A(\overline{K}) \overset{2}{\longrightarrow} A(\overline{K}) \to 0.
$$ 
Let $A = J(C)$ be the Jacobian of the genus 2 curve $C\!: y^2 = f(x)$ over $K$.  As is well-known \cite[Ch.\! 2]{CF}, the effective divisor $\gO$ of degree 2 on $C$ lying over $x = \infty$ represents the canonical divisor class and every $K$-rational divisor class of degree 0 is represented by a unique divisor of the form $\gA - \gO$, with $\gA$ effective of degree 2 and defined over $K$.   For our applications,  assume that $\deg f = 6$ and factor 
$$
f(x) = c(x-r_1) \cdots (x-r_6) 
$$
over the splitting field $F$ of $f$.   Let $\cH_F$ be the quotient of $\oplus_{j=1}^6  \, F^{\times}/F^{\times 2}$ by the image of $F^\times$ on the diagonal.  If $\gA = P_1+P_2$, where the points $P_i = (x_i,y_i)$ are affine and $x_i$ is not a root of $f$, then the $j$-th coordinate of the $x-T$ homomorphism $\partial_F\!: A(F) \to \cH_F$ is induced by ${\gA}-{\gO} \leadsto (x_1- r_j)(x_2-r_j)$.  However, if say $x_1 = r_j$, then a non-zero entry is obtained by replacing $x_1-r_j$ by $f'(r_j)$, since
$$
\frac{y^2}{x-r_j} = \frac{f(x)}{x-r_j} = c \prod_{\begin{smallmatrix} i= 1, \dots, 6 \\ i \ne j \end{smallmatrix}} (x-r_i).
$$ 
The kernel of $\partial_F$ is $2A[F]$.   If $\partial_F(P)$ is represented by $(c_1, \dots, c_6)$ in $\oplus \, \F^{\times}$, then the product $c_1 \cdots c_6$ is in $F^{\times 2}$.

\begin{lem} \label{x-T}
Let $P$ in $A[2]$ be represented by the divisor $(r_5,0)+(r_6,0) - \gO$ and let $q(x) = (x-r_5)(x-r_6)$.  Then $\partial_F(P)$ is represented by 
$$
\left(q(r_1), \hspace{4 pt} q(r_2), \hspace{4 pt} q(r_3), \hspace{4 pt} q(r_4), \hspace{4 pt}  (r_6-r_5)f'(r_5), \hspace{4 pt} (r_5-r_6)f'(r_6)\right).
$$
The elementary $2$-extension $F(\frac{1}{2}P)/F$ is generated by the square roots of
$$
 \frac{c(r_5-r_1)(r_6-r_2)}{(r_6-r_3)(r_6-r_4)}, \, \frac{c(r_5-r_2)(r_6-r_1)}{(r_6-r_3)(r_6-r_4)}, \, \frac{c(r_5-r_3)(r_6-r_1)}{(r_6-r_2)(r_6-r_4)},  \, \frac{c(r_5-r_4)(r_6-r_1)}{(r_6-r_2)(r_6-r_3)}.
$$
\end{lem}

\begin{proof}
We find $\partial_F(P)$ by the definition reviewed above.  Since the diagonal image of $F^\times$ is trivial in $\cH_K$, we divide by the last coordinate.  Then the product of the first five coordinates is a square, so adjoining the square roots of the first four coordinates is necessary and sufficient to obtain $F(\frac{1}{2}P)$.   
\end{proof}

Now let $k$ be a finite field of characteristic 2, let $K$ be the field of fractions of $\W = \W(k)$ and let $A = J(C)$ as in Proposition \ref{FindParams}.    Replace $x$ by a scalar multiple to arrange that $a_5$ be one of (at most 5) fixed representatives for $\W^{\times}/\W^{\times 5}$ and  translate to make the coefficient of $x^4$ zero. 
Then $C$ has  a  model of the form $y^2 = x \Phi(x)$, where $\Phi(x)$ is the quintic
\begin{equation} \label{PhiEq}
\Phi(x) =  x^5(1 + 4g(1/x)) = (4a_0+1)x^5+4a_1x^4+4a_2x^3+4a_3 x^2+4a_5.
\end{equation}
Let $r_6 = 0$ and let $r_1, \dots, r_5$ be the roots of $\Phi$.

\begin{lem} \label{FindF}
The $2$-division field $F$ of $A$ is $K(\zeta, \pi)$, where $\zeta$ is a primitive fifth root of unity and $\pi$ is a root of $a_5^2 z^5 - 2$.  The roots of $\Phi$ have the form 
$$
r_j= -a_5 \zeta^j \pi^2 + 0(2\pi), \quad j = 1, \dots, 5.
$$  
\end{lem}

\begin{proof}  
Let $L = K(\zeta, \pi)$ and $d(z) = \frac{1}{\pi^{10}} \Phi(\pi^2 z) \equiv z^5 + a_5^5 \pmod{ \pi^4 \cO_L[z] }$.  By Hensel's Lemma, the approximate root $- a_5\zeta^j$ leads to a root of $d(z)$ of the form $z_j = - a_5\zeta^j +0(\pi^4)$ and so $r_j = z_j \pi^2$ is in $L$.  Conversely, any root $r$ of $\Phi(x)$ satisfies $r^5 \equiv -4a_5 \pmod{4r^2}$ and $\ord_2(r) = 2/5$, so $\pi_0 = 2/r^2$ leads to a root of $a_5^2 z^5 - 2$ by Hensel.  By Proposition \ref{FieldForE}, $\zeta_5$ is in $F$, so $L=F$.  
\end{proof}

Let $\tilde{C}\!: \, y^2 + y = \tilde{g}(x)$ be the deformation of $C$ with $\tilde{g}$ in equation \eqref{tildeg} and 
\begin{equation} \label{perturb}
b_3 = a_3 +4 \epsilon_3,  \quad  b_j = a_j + 2\epsilon_j  
\text{ for } j \in \{1,2\}, \quad b_0 = a_0 + \epsilon_0.
\end{equation}
Note that the leading coefficent $a_5$ of $\tilde{g}$ agrees with that of $g$, since it has been chosen from the discrete set of representatives for $\W^{\times}/\W^{\times 5}$. By the Lemma, $A$ and $\tilde{A}$ have the same 2-division field $F$.  Let 
$
\tilde{\Phi} =  (4b_0+1)z^5+4b_1z^4+4b_2z^3+4b_3 z^2+4a_5.
$  

\begin{lem} \label{tilder}
For each root $r$ of $\Phi$, there is a root $\tilde{r}$ of $\tilde{\Phi}$  and an $\epsilon_0$ in $\cO_F$ such that 
$
\tilde{r}/r \equiv 1+ 4\epsilon_0  \pmod{4 \pi}.
$
If $r' \ne r$ is another root of $\Phi$, then 
$$
(\tilde{r}'-\tilde{r})/(r'-r)  \equiv 1+ 4\epsilon_0  \pmod{4 \pi}.
$$  
Let $P$ in $A[2]$ be represented by $(r,0)+(r',0) - \gO$ and let $\tilde{P}$ be the corresponding point in $\tilde{A}[2]$.  Then $F(\frac{1}{2}\tilde{P}) = F(\frac{1}{2}P)$.
\end{lem}

\begin{proof}
To describe the roots of $\tilde{\Phi}$ as multiples of the roots of $\Phi$, let
$$
h(z) = \frac{1}{r^5} \tilde{\Phi}(rz) = (4b_0+1) z^5 + \frac{4b_1}{r} z^4 + \frac{4b_2}{r^2} z^3 + \frac{4b_3}{r^3} z^2 + \frac{4a_5}{r^5}.
$$ 
Then $h(z)$ is in $\cO_F[z]$ and $h(z) = \frac{1}{r^5} \Phi(rz) + \delta(z)$, where
$$
\delta(z) = 4\epsilon_0 z^5 + \frac{8\epsilon_1}{r} z^4 + \frac{8\epsilon_2}{r^2} z^3 + \frac{16\epsilon_3}{r^3} z^2 \, \equiv \, 4\epsilon_0 z^5 \pmod{4\pi \, \cO_F[z]}. 
$$
Since $h'(1) \equiv 1 \pmod{2\pi}$ is a unit in $\cO_F$ and $h(1) = \delta(1) \equiv 0 \pmod{4}$, Hensel's Lemma implies that there is a root $s$ of $h(z)$ satisfying $s \equiv 1 \bmod{4}$.  By quadratic convergence of Hensel iteration we find that  
$$
s \equiv 1 - \frac{h(1)}{h'(1)} \equiv   1 + 4\epsilon_0  \pmod{4\pi}.  \hspace{20 pt} 
$$
Finally, modulo $4 \pi$, the ratio of each Kummer generator in Lemma \ref{x-T} to the corresponding Kummer generator arising from $\tilde{P}$ is $(1+4 \epsilon_0)^4 \equiv 1 \pmod{4 \pi}$ and therefore a square in $F^\times$.  Hence $F(\frac{1}{2}P) = F(\frac{1}{2}\tilde{P})$.
\end{proof}

We briefly indicate how the Kummer groups of curves and  Honda systems are related.   Let $\cE$ be the group scheme belonging to the finite Honda system $\gE=\gE_1$ in Notation \ref{Ebasis} with $k = \F_2$.   By Proposition \ref{FieldForE}(iii), each non-trivial point $P_a$ in $E$ corresponds to a root $a$ of $x^{15}+8$.  Let $\cV$ be the group scheme belonging to the finite Honda system $\gM$ of exponent 4 in \S \ref{expp2} with Honda parameters $\bfs = [\lambda;s_1,s_2,s_3,s_5]$ as in Definition \ref{HParam}.  Thus $\lambda \equiv \pm 1 \bmod{4}$ and $s_j$ is in $\F_2$.  The extension $L_a = F(\frac{1}{2}P_a)$ of $F$ generated by $Q$ in $V$ such that $2Q = P_a$ does not depend on the choice of $Q$. In addition, $L_a/F$ is an elementary 2-extension whose degree divides 16.

As in the proof of Lemma \ref{lem:L01}, we found a polynomial $f_{\bfs,a}(x)$ in $F[x]$ with splitting field $L_a$.  For each $a$ and each of the 32 Honda parameters $\bfs$, we used Magma to obtain the {\em Kummer group} for $L_a$, i.e.\! the subgroup of $F^{\times}/F^{\times 2}$ whose square roots generate $L_a$ over $F$.  

Next, we match $P_a$ coming from the finite Honda system to a point of order 2 on the Jacobian $A$ of the curve $C\!: \, y^2 = x \Phi(x)$ with $\Phi$ given by \eqref{PhiEq}.  The field $F = \Q_2(\pi,\zeta)$ has an explicit construction, where $\zeta$ is a primitive fifth root of unity and $\pi$ is a prime element satisfying $\pi^5 = 2$.   Then there is a Frobenius $\tau$ in $\Gal(F/\Q_2)$ fixing $\pi$ and a generator $\sigma$ for the inertia group in $\Gal(F/\Q_2)$ such that $\sigma(\pi) = \zeta \pi$.   Since $-\pi$ is the unique root of $x^{15} + 8$ fixed by $\tau$ and the divisor $(r_5,0)+(0,0) - \gO$ represents the unique point $T$ in $A[2]$ fixed by $\tau$, we find that $P_{-\pi}$ and $T$ correspond.  

By Lemma \ref{tilder}, it suffices to determine the Honda parameters for the Jacobians of a finite set of curves.  For each of these curves, we compare the Kummer group in Lemma \ref{x-T} to the Kummer groups obtained by factoring $f_{\bfs,-\pi}$ and thereby determine the appropriate value of $\bfs$. 

\section{Global Applications}  \label{global}  \numberwithin{equation}{section}
Throughout this section, $F$ is the Galois closure of a favorable quintic field of discriminant $\pm 16N$, with $N$ prime, as in Definition \ref{fav}.

\begin{Rem}   \label{favS}
If $A$ is a favorable abelian surface of conductor $N$, then $\Q(A[4])$ is a favorable $\gS$-field containing $F = \Q(A[2]),$ as  in Definition \ref{Sfield}.  Item (i) is standard, (ii) follows from Proposition \ref{CondBd} and (iii)  is in Proposition \ref{BigGal} and its proof.
\end{Rem}

Since $\vert \gS \vert= 120\cdot 2^{11}$, a realistic test for the existence of a favorable $\gS$-field $L$ is desirable. We describe a subfield $K'$ of $L$ whose Galois closure is $L$, i.e. a {\em stem field} for $L$.  First, let $K$ be the subfield of $F$ fixed by $\ov{H} = \Sym \{1,2,3\}$, so that $K$ is obtained by symmetrizing $r_5-r_4$, where $r_1, \dots, r_5$ are the roots of a favorable quintic polynomial with splitting field $F$.  The transposition $s = (45)$ commutes with $\Gal(F/K)$ and so  induces an automorphism of $K$.  A standard double coset computation shows that there is a unique prime $\gn_s$ over $N$ in $K$ fixed by $s$.  Moreover, there is a unique prime $\gp$ over 2 in $K$.

\begin{prop} \label{STest}
Let $L$ be a favorable $\gS$-field.  Then $L$ admits four stem fields $K'$ such that $K'/K$ is quadratic with ray class modulus $\gp^6 \gn_s \infty$.
\end{prop}

\begin{proof}
Let $H = \pi^{-1}(\ov{H})$ in $\gS$.   There are 31 subgroups $H'$ of index 2 in $H$.  The action of $\gS$ on left cosets of $H'$ is faithful for 12 of them.   For those, the fixed field $K'$ of $H'$ is a stem field for $L$.  Let $\sigma$ be a transvection in $\gS$ satisfying $\pi(\sigma) = (45)$.  The conjugates of $\sigma$ by representatives $t$ in $\gS$ for the double cosets $Ht\lr{\sigma}$ generate the inertia groups at primes over $N$.  Intersecting these inertia groups with $H$ and $H',$ respectively, we find that for four choices of $K'/K$, the only ramification over $N$ occurs over $\gn_s$.  Definition \ref{Sfield}(ii) provides the conductor bound over 2. 
\end{proof}

\begin{cor} \label{STestCor}
Let $F$ be the Galois closure of a favorable quintic field of discriminant $\pm16N$ and let $K$ be the subfield of $F$ fixed by $\Sym\{1,2,3\}$.  If no quadratic extension $K'/K$ of modulus $\gp^6 \gn_s \infty$ exists such that the Galois closure of $K'/\Q$ has Galois group $\gS$,  then no abelian surface of conductor $N$ with 2-division field $F$ exists.
\end{cor}

\begin{Rem}  \label{mN}
The non-existence results in \cite[Table 3]{BK3} follow from the Corollary.  Under GRH, it was used to test all favorable quintic fields in the Bordeaux tables \cite{BT}.  For 311 primes $N$ in $\{1277,1597,2557, \dots, 310547,312413\}$, we found that no favorable abelian surface of conductor $N.$ 

Corollary \ref{STestCor} admits a stronger conclusion.  Suppose that $B_{/\Q}$ is a semistable abelian surface of conductor $qN$, with $q$ odd and prime to $N$, that $\Q(B[2])= F$ and that the primes over $q$ are unramified in the $4$-division field $L$ of $B$.  Then $L$ is a favorable $\gS$-field.  If no favorable $\gS$-field containing $F$ exists, then all such $B$ are ruled out, not just those with $q = 1$.

In fact, such $B$  exist with $q > 1$.  In Table \ref{Nmtable}, $[a_0, \dots, a_6]$ denotes a polynomial $f(x) = a_0 + \dots + a_5 x^5$ defining a  favorable quintic field.  The Igusa discriminant of a minimal model for the curve $y^2 = f(x)$ is $\Delta = q^4 N$ and its Jacobian $B$ is semistable of conductor $qN$.  Let $\cB$ be the N\'eron model of $B$ over the strict Henselization $\cO$ of $\Z_q$ and let $k \simeq \ov{\F}_q$ be the residue field of $\cO$.  Let $\cB^0_k$ be the connected component of the identity in the special fiber $\cB_k$.  Then the group of connected components $\Phi = \cB_k/\cB^0_k$ is cyclic of order $4 = \ord_q(\Delta)$ by \cite[Prop.\! 2(ii)]{L}.  Moreover, $\cB^0_k$ is the extension  of an elliptic curve by a torus of dimension 1, so $\cB^0_k[4]  \simeq (\Z/4\Z)^3$.  Since the kernel of reduction is divisible by 4, all these points lift to $\cB(\cO)[4]$.  Hence $q$ is unramified in $\Q(B[4])$.

The examples in Table \ref{Nmtable} occurred in an old collection of abelian surfaces created by the senior author, rather than a dedicated search. They  illustrate why the criterion in Corollary \ref{STestCor} often is inadequate to decide the non-existence of a favorable abelian surface of conductor $N$ in \cite[Table 3]{BK3}.  Conversely, given a favorable $\gS$-field $L$, is there always a semistable abelian surface $B$ of conductor $qN$ with $\Q(B[4)) = L$ for some odd integer $q$?

\begin{table}[h]\begin{tabular}{|c|c|c|}
\hline
$N$&$q$&$f(x)$\\
\hline 
1061&3&[ 1, 8, 28, 56, 64, 36 ]\\
2069&31&[ 25, 344, 1888, 5168, 7056, 3844 ]\\
2269&3&[ 1, 12, 56, 124, 124, 36 ]\\
2909&3&[ 1, 0, -12, -4, 32, 36 ]\\
3989&11&[ -11, -112, -404, -532, 68, 484 ]\\
5381&5&[ 1, 8, -8, 0, 0, 4 ]\\
7013&5&[ 1, 0, 12, 16, 20, 4 ]\\
7877&3&[ 1, -8, 16, 12, -56, 36 ]\\
8581&11&[ -3, 52, -324, 876, -1068, 484 ]\\
\hline
\end{tabular}
\smallskip
\begin{caption}{
$B$ has conductor $qN$ and $\Q(B[4])$ is a favorable $\gS$-field} \label{Nmtable}\end{caption}
\end{table}
\end{Rem}

\noindent{\bf Errata}.  We  fix a misquote in \cite{BK2}. The orthogonal group $O^+_4(\F_2) \simeq \cS_3 \wr \cS_2$ is not generated by transvections and  does not belong  in Propositions 2.4 and 2.5.  In Propositions 2.8 and 2.11, $V$ should be semistable, so that Fontaine's bound on ramification at $\ell$ applies.

\end{document}